\documentclass[12pt]{article}
\usepackage[margin=1in]{geometry}
\usepackage[T1]{fontenc}
\usepackage{amsmath,amssymb,amsthm}
\usepackage{mathpazo}
\usepackage{graphicx}

\usepackage{setspace}
\usepackage{microtype}
\usepackage{needspace}
\usepackage[authoryear,round]{natbib}
\usepackage[hidelinks]{hyperref}
\newtheorem{theorem}{Theorem}
\newtheorem{proposition}[theorem]{Proposition}
\newtheorem{lemma}[theorem]{Lemma}
\newtheorem{corollary}[theorem]{Corollary}
\theoremstyle{remark}\newtheorem{remark}[theorem]{Remark}
\newcommand{\E}{\mathbb E}
\newcommand{\Pp}{\mathbb P}
\newcommand{\R}{\mathbb R}
\newcommand{\C}{\mathbb C}
\newcommand{\ii}{\mathrm i}
\DeclareMathOperator{\Ree}{Re}

\usepackage{xcolor}
\definecolor{darkgreen}{rgb}{.2,0.4,0.3}

\title{Analyticity of $\operatorname{Sine}_\beta$ pair correlations}
\author{Yun Li,  Yahui Qu, Benedek Valk\'o, B\'alint Vir\'ag}
\hypersetup{pdfauthor={Yun Li,  Yahui Qu, Benedek Valk\'o, B\'alint Vir\'ag},
  pdftitle={Analyticity of Sine-β pair correlations}}
\date{}
\begin{document}
\maketitle
\begin{abstract}
We prove that the pair correlation function of $\operatorname{Sine}_\beta$
extends holomorphically in $\beta$ to $\Ree\beta>0$ at every fixed nonzero
separation. Our proof combines a random walk representation of the
stochastic sine phase moments, stretched-exponential Fourier decay
obtained by coupling, and weighted energy estimates for complex $\beta$.
\end{abstract}

\section{Introduction}
Let $\Xi_\beta$ be the $\operatorname{Sine}_\beta$ point process, normalized
to have intensity $1/(2\pi)$. Its pair correlation function $\rho_\beta$
is defined by
\[
 \E\sum_{x,y\in\Xi_\beta}^{\ne}F(x,y)
 =\int_{\R^2}F(x,y)\rho_\beta(x,y)\,dx\,dy
\]
for nonnegative compactly supported Borel functions $F$, where the sum
runs over ordered pairs of distinct points.

\begin{theorem}[Analyticity of the pair correlation]\label{thm:main} For every $\beta>0$, the pair correlation function $\rho_\beta$
has a continuous version off the diagonal. For this version and
each $\lambda>0$, the function
\[
 \beta\longmapsto\rho_\beta(0,\lambda)
\]
extends holomorphically from $(0,\infty)$ to $\{\Ree\beta>0\}$.
The extension is jointly continuous in $(\beta,\lambda)$ on
$\{\Ree\beta>0\}\times(0,\infty)$.
\end{theorem}

\begin{figure}[t]
\centering
\includegraphics[width=\textwidth]{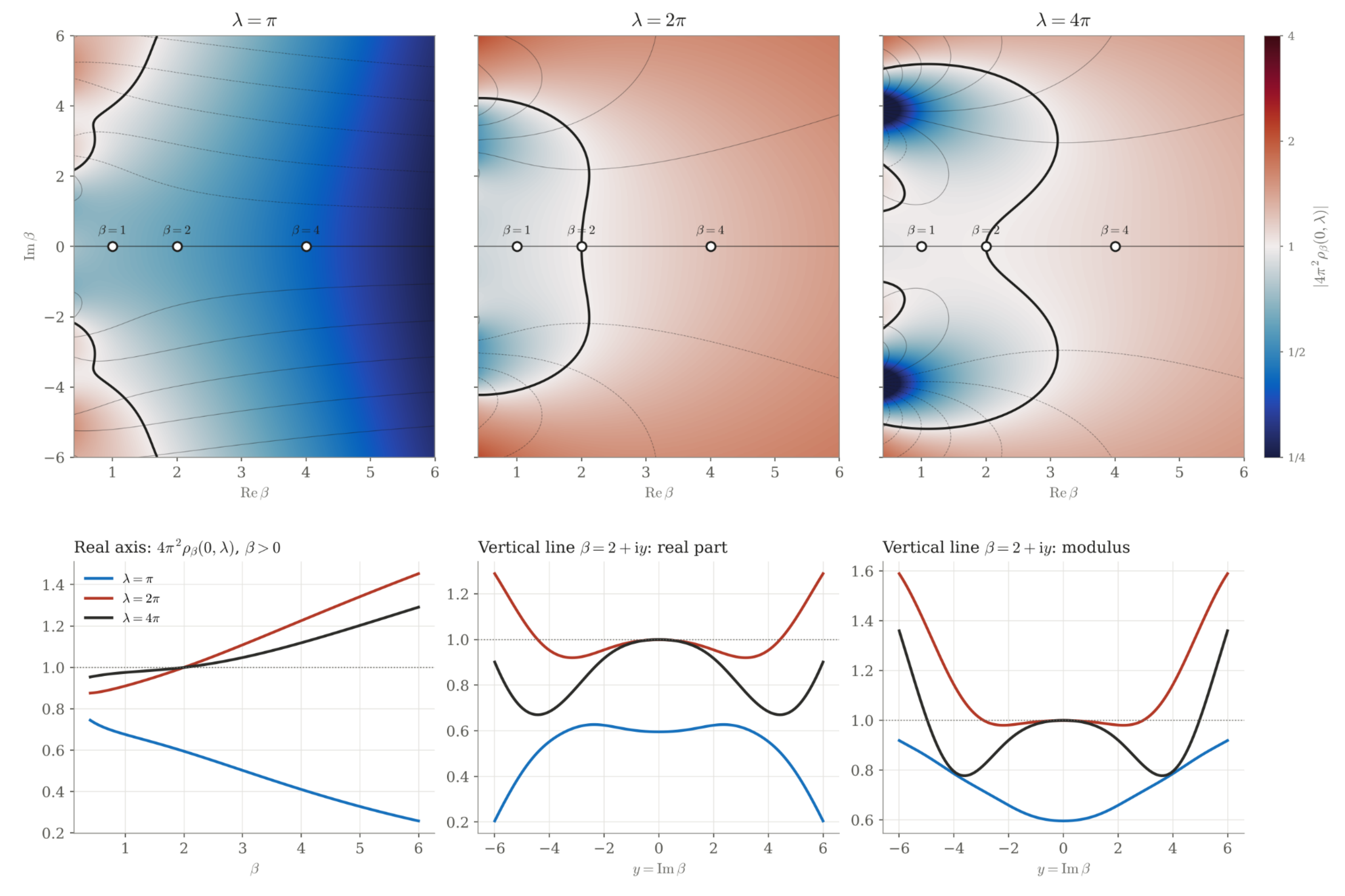}
\caption{The holomorphic continuation of $4\pi^2\rho_\beta(0,\lambda)$, computed numerically
from \eqref{eq:momentode}. The black curves on the contour plots are where $4\pi^2\left|\rho_\beta(0,\lambda)\right|=1$.}
\label{fig:complexpair}
\end{figure}

Figure~\ref{fig:complexpair} shows the continuation numerically.

Dependence on inverse temperature has long been studied through exact
formulas and perturbation theory. For $\beta=1,2,4$, the correlation functions of
$\operatorname{Sine}_\beta$ have classical explicit formulas:
determinantal for $\beta=2$ and Pfaffian for $\beta=1,4$
\citep[Chapters~5--7]{Forrester10}. \cite{Forrester92,Forrester94} obtained correlation formulas
at even integer $\beta$, and \cite{Ha95} extended
these to rational coupling.  \cite{Okounkov97} subsequently obtained Selberg-type integral formulas for
$n$-point correlations at rational $\beta$. \cite{FJM01}
computed first-order perturbations around $\beta=2$ and $4$, following earlier work on the two-dimensional
plasma \citep{Jancovici81}. More recently, \cite{QV25} obtained a
stochastic representation for all $\beta>0$ and proved continuity in
$\beta$. Theorem~\ref{thm:main} strengthens continuity to
holomorphic dependence throughout the right half-plane.

Translation and reflection invariance reduce the pair correlation to
$\rho_\beta(0,\lambda)$, $\lambda>0$.
Our proof uses the stochastic phase representation of
$\operatorname{Sine}_\beta$ \citep{KS09}. We express the pair
correlation as the phase density at zero divided by $2\pi$, and use a
random walk representation of the phase moments to prove analyticity.

Let $Z=B_1+\ii B_2$, where $B_1,B_2$ are
independent two-sided standard real Brownian motions. For real $\beta>0$
and $\lambda\in\R$, let $\alpha_{\lambda,\beta}$ solve the phase SDE
\begin{equation}\label{eq:phase}
 \begin{gathered}
 d\alpha_{\lambda,\beta}(t)
 =\frac{\beta\lambda}{4}e^{\beta t/4}\,dt+
 \Ree\!\left[(e^{-\ii\alpha_{\lambda,\beta}(t)}-1)\,dZ_t\right],\qquad t\le0,\\
 \lim_{t\to-\infty}\alpha_{\lambda,\beta}(t)=0.
 \end{gathered}
\end{equation}
Equation~\eqref{eq:phase} is the Killip--Stoiciu equation in logarithmic
time, as considered in \citet[Corollary~51]{VV22}.
All spectral parameters are driven by the same $Z$. We often write
$\alpha_\lambda$ when $\beta$ is fixed or clear from context.
If $U$ is uniform on $[0,2\pi)$ and independent of the phase field, then
\begin{equation}\label{eq:pointprocess}
 \Xi_\beta=\{\lambda\in\R:
                   \alpha_{\lambda,\beta}(0)=U\bmod2\pi\}
\end{equation}
has the $\operatorname{Sine}_\beta$ distribution
\citep[Definition~1.10 and Section~4]{KS09}.
Write $\mathbb T=\R/(2\pi\mathbb Z)$ and let $p_{\beta,\lambda}$
denote the density of $\alpha_{\lambda,\beta}(0)\bmod 2\pi$ with respect
to $d\theta$ on $\mathbb T$, for $\lambda>0$.
For $k\in\mathbb Z$, define its Fourier coefficients
\[
 q_k(\beta,\lambda)=\E e^{\ii k\alpha_{\lambda,\beta}(0)}.
\]

\begin{theorem}[Pair correlation and phase density]\label{prop:pairphase}
For real $\beta,\lambda>0$, the pair correlation off the diagonal and
the phase density have continuous versions, and
\begin{equation}\label{eq:pairfourier}
 \rho_\beta(0,\lambda)
 =\frac{1}{2\pi}\,p_{\beta,\lambda}(0)
 =\frac{1}{4\pi^2}\sum_{k\in\mathbb Z}q_k(\beta,\lambda),
\end{equation}
with absolute convergence.
\end{theorem}
The existence of correlation functions of all orders for
$\operatorname{Sine}_\beta$, $\beta>0$, was proved by \citet{QV25}.
More recently, \citet{AN26} represented these functions in terms of
the Hua--Pickrell stochastic zeta function introduced by \citet{LV_HP}, and established their
continuity and boundedness in the spatial variables.
For higher correlations, we obtain the following representation
at the level of correlation measures.

\begin{proposition}[Phase representation of correlation measures]\label{prop:anglelimit}
Let $\beta>0$, $k\ge2$, and let $I_0,\ldots,I_{k-1}$ be pairwise disjoint bounded
intervals. Write $\alpha_\lambda=\alpha_{\lambda,\beta}(0)$. Then
\begin{equation}\label{eq:anglelimit}
 \E\prod_{j=0}^{k-1}\Xi_\beta(I_j)
 =\frac1{2\pi}\lim_{\varepsilon\downarrow0}
 \int_{I_0\times\cdots\times I_{k-1}}
 \frac{\Pp\!\left({\max\limits_{1\le j\le k-1}}
             \operatorname{dist}(\alpha_{\lambda_j-\lambda_0},2\pi\mathbb Z)
                         <\varepsilon\right)}
      {(2\varepsilon)^{k-1}}\,d\lambda_0\cdots d\lambda_{k-1} .
\end{equation}
\end{proposition}
The proof is in Section~\ref{sec:pair}. When $k=2$, continuity of the phase
density turns this limit into \eqref{eq:pairfourier}.

It\^o's formula and the scaling invariance of the phase SDE
\eqref{eq:phase} yield
\begin{equation}\label{eq:momentode}
 \lambda\partial_\lambda q_k
 =\ii k\lambda q_k+
 \frac{2k^2}{\beta}(q_{k+1}+q_{k-1}-2q_k),
 \qquad q_0=1,
\end{equation}
where all coefficients are evaluated at $(\beta,\lambda)$.
The derivation is given in Section~\ref{sec:walk}.
To extend \eqref{eq:pairfourier} to complex $\beta$, we first continue
each coefficient and then establish locally uniform summability. The key
tool is a representation by a random walk whose law does not depend on
$\beta$.

\begin{proposition}[Random walk representation]\label{prop:walk}
Let $(S_t)_{t\ge0}$ be the continuous-time nearest-neighbor random walk on
$\mathbb Z_{\ge0}$, absorbed at zero. From $k\ge1$, its jump rate to
each of $k-1$ and $k+1$ is $k^2/2$.
Let $\E_k$ denote expectation for the walk started at $k$.
For real $\beta>0$, $\lambda>0$,
and $k\ge1$,
\begin{equation}\label{eq:area}
 q_k(\beta,\lambda)=\E_k e^{\ii\lambda A_\beta},\qquad
 A_\beta=\frac\beta4\int_0^\infty e^{-\beta t/4}S_t\,dt.
\end{equation}
\end{proposition}

Section~\ref{sec:walk} proves this formula by applying the countable-state
Feynman--Kac formula on $[-T,0]$ and then sending $T\to\infty$. The embedded
simple random walk is absorbed after finitely many jumps almost surely,
so $S$ is nonexplosive and the integral in \eqref{eq:area} is almost surely
finite. Write $\tau_j=\inf\{t\ge0:S_t=j\}$. For $1\le k\le n$,
the finite-space cutoff is
\begin{equation}\label{eq:cutoff}
 q_{k,n}(\beta,\lambda)
 =\E_k\!\left[e^{\ii\lambda A_\beta}
             \mathbf1_{\{\tau_0<\tau_{n+1}\}}\right].
\end{equation}
Gambler's ruin gives the uniform real-parameter estimate
\begin{equation}\label{eq:ruin}
 |q_k(\beta,\lambda)-q_{k,n}(\beta,\lambda)|
 \le\frac{k}{n+1},\qquad \beta>0,\quad\lambda>0.
\end{equation}
The cutoffs are holomorphic on $\Ree\beta>0$. A weighted energy estimate
provides a bound independent of $n$, allowing Vitali's theorem to remove
the cutoff. The details are in Section~\ref{sec:holomorphic}.

\begin{proposition}[Holomorphic Fourier coefficients]\label{prop:holomorphic}
For every fixed $\lambda>0$ and $k\in\mathbb Z$, $q_k(\beta,\lambda)$ has a unique
holomorphic extension to $\Ree\beta>0$. These extensions satisfy
\begin{equation}\label{eq:complexbound}
 \sum_{k=1}^{\infty}\frac{|q_k(\beta,\lambda)|^2}{k^2}
 \le 4\frac{|\beta|^2}{(\Ree\beta)^2},
 \qquad \lambda>0,\quad \Ree\beta>0.
\end{equation}
Here $q_0(\beta,\lambda)=1$ and
$q_{-k}(\beta,\lambda)=\overline{q_k(\bar\beta,\lambda)}$ for $k\ge1$.
\end{proposition}

To sum these holomorphic extensions in \eqref{eq:pairfourier}, we need
locally uniform bounds that are summable in $k$.
For real $\beta>0$, a coupling argument for the random walk in
Proposition~\ref{prop:walk} gives stretched-exponential decay
(Section~\ref{sec:real}). Combining this decay with the energy bound
\eqref{eq:complexbound}, we obtain the same type of estimate for
complex $\beta$ with $\Ree\beta>0$ (Section~\ref{sec:transfer}):

\begin{proposition}[Stretched-exponential Fourier decay for complex $\beta$]\label{prop:rapid}
For every compact $K\subset\{\Ree\beta>0\}$ and compact interval
$L\subset(0,\infty)$, there are constants $a,b>0$ such that
\[
 \sup_{\substack{\beta\in K\\\lambda\in L}}|q_k(\beta,\lambda)|
 \le b e^{-a|k|^{1/3}},\qquad k\in\mathbb Z.
\]
\end{proposition}

Thus the series in \eqref{eq:pairfourier} converges locally uniformly for
$\Ree\beta>0$ and gives the extension in Theorem~\ref{thm:main}.

Section~\ref{sec:preliminaries} records the phase representation and
the probabilistic and analytic preliminaries. Sections~\ref{sec:walk}
and~\ref{sec:real} establish the random walk representation and real
Fourier decay. Section~\ref{sec:pair} identifies the pair correlation.
Sections~\ref{sec:holomorphic} and~\ref{sec:transfer} complete the
analytic continuation argument.
Section~\ref{sec:beta2decay} gives explicit formulas for the Fourier
moments and phase density at $\beta=2$, recovering the sine-kernel pair
correlation. The resulting moment asymptotics show that the decay exponent
$1/3$ in Proposition~\ref{prop:rapid} cannot be increased beyond $1/2$
in general.

\section{Preliminaries}\label{sec:preliminaries}
\subsection{The phase and its boundary angle}
\begin{proposition}[The stochastic sine phase]\label{prop:shooting}
For real $\beta>0$, the phase SDE \eqref{eq:phase} has a pathwise unique
solution for the prescribed driving Brownian motion.
Almost surely $\lambda\mapsto\alpha_\lambda(0)$ is a real analytic increasing
bijection of $\R$, with $\alpha_0(0)=0$ and
$\partial_\lambda\alpha_\lambda(0)>0$.
Moreover, for $\lambda\in\R$ and $t\le0$,
\begin{equation}\label{eq:phasemean}
 \alpha_\lambda(t)\stackrel d=\alpha_{\lambda e^{\beta t/4}}(0),
 \qquad \E\alpha_\lambda(t)=\lambda e^{\beta t/4}.
\end{equation}
Let $U$ be uniform on $[0,2\pi)$ and independent of $Z$, and define
$\Xi_\beta$ by \eqref{eq:pointprocess}.
Then $\Xi_\beta$ has the $\operatorname{Sine}_\beta$ distribution.
The corresponding normalized spectral weights
\[
 w_\lambda=\frac{1}{\partial_\lambda\alpha_\lambda(0)},
 \qquad \lambda\in\Xi_\beta,
\]
are mutually independent gamma variables of shape $\beta/2$ and mean one,
and are independent of $\Xi_\beta$.
\end{proposition}
These weights are half the right spectral weights in
\citet[Proposition~3 and Lemma~17]{VV23}.
The phase and shooting representations follow from
\citet[Corollary~51]{VV22}; see also \citet[Definition~16 and Lemma~17]{VV23}.
Time translation in \eqref{eq:phase} gives the scaling in
\eqref{eq:phasemean}. The expectation identity follows from
\citet[Proposition~4.5]{KS09} after the time change $s=e^{\beta t/4}$.

The phase representation was introduced by \citet{KS09}; its equivalence
with the Gaussian bulk limit of \citet{VV09} was proved by
\citet[Corollary~1.7]{Nakano14} and \citet[Theorem~28]{VV17}.

\begin{lemma}[Phase and counting moments]\label{lem:phasebounds}
For each $0<L<\infty$ there is $\eta_L>0$ such that
\begin{equation}\label{eq:phasebounds}
 \sup_{\beta>0,\,0\le\lambda\le L,\,t\le0}
       \E e^{\eta_L\alpha_\lambda(t)}<\infty.
\end{equation}
Thus phase moments of every positive order are locally bounded in $\lambda$,
uniformly in $\beta>0$ and $t\le0$. Counts in bounded intervals have moments
of every positive order as well.
\end{lemma}
\begin{proof}
Put $N_\lambda=\Xi_\beta((0,\lambda])$. The counting bound in
\citet[Proposition~9(viii)]{VV09}, applied with total
drift $\lambda$, gives
\[
 \Pp(N_\lambda\ge ak)\le2\left(\frac{\lambda}{2\pi a}\right)^k,
 \qquad a,k\in\mathbb N.
\]
Choose an integer $a>L/(2\pi)$. This gives a positive exponential moment
for $N_\lambda$ uniformly for $\lambda\le L$ and $\beta>0$.
By \eqref{eq:pointprocess}, almost surely,
\begin{equation}\label{eq:rounding}
 N_\lambda=\left\lfloor\frac{\alpha_\lambda(0)+2\pi-U}{2\pi}\right\rfloor.
\end{equation}
Consequently, $|2\pi N_\lambda-\alpha_\lambda(0)|\le2\pi$,
so the phase has the same moment property. The scaling identity in
\eqref{eq:phasemean} gives the bound for $t\le0$. Stationarity gives
the count assertion for any bounded interval.
\end{proof}

\subsection{Probabilistic and analytic tools}\label{sec:probtools}

\begin{lemma}[Feynman--Kac]\label{lem:fk}
Let $S$ be a nonexplosive continuous-time Markov chain on a countable
state space $E$, with finite total jump rates and generator $\mathcal G$.
Fix $T>0$, and let $v:[0,T]\times E\to\C$ be continuous in time,
with $\Ree v(t,k)\le C$ for some $C\ge0$.
Every bounded solution $u:[0,T]\times E\to\C$ of
\begin{equation}\label{eq:fk-system}
 \partial_tu(t,\cdot)=\mathcal G u(t,\cdot)+v(t,\cdot)u(t,\cdot),
 \qquad u(0,\cdot)=g,
\end{equation}
with $u(\cdot,k)\in C^1([0,T])$ for each $k\in E$, satisfies
\begin{equation}\label{eq:fk-formula}
 u(t,k)=\E_k\!\left[
 \exp\!\left(\int_0^t v(t-s,S_s)\,ds\right)g(S_t)\right],
 \qquad 0\le t\le T.
\end{equation}
In particular, there is at most one bounded solution.
\end{lemma}
\begin{proof}
Fix $t\le T$, and set
\[
 M_r=\exp\!\left(\int_0^r v(t-s,S_s)\,ds\right)u(t-r,S_r),
 \qquad 0\le r\le t,
\]
and let $\sigma_n$ be the exit times from increasing finite subsets
exhausting $E$. The jump-process product rule shows that
$M_{r\wedge\sigma_n}$ is a martingale, since its drift before
$\sigma_n$ is proportional to $-\partial_tu+\mathcal Gu+vu=0$.
The finite jump rates and boundedness of $u$ justify the generator sums.
Moreover,
\[
 |M_{r\wedge\sigma_n}|\le e^{CT}\|u\|_\infty.
\]
Nonexplosion gives $\sigma_n\to\infty$ almost surely, so bounded
convergence yields $\E_k M_t=M_0$, proving \eqref{eq:fk-formula}.
See also \citet{FP99} for general Markov-process formulations.
\end{proof}

\begin{lemma}[Hardy's inequality in difference form]\label{lem:hardy}
For $z_1,\ldots,z_n\in\C$, with $z_0=z_{n+1}=0$,
\begin{equation}\label{eq:hardy}
 \sum_{k=1}^n\frac{|z_k|^2}{k^2}
 \le4\sum_{k=0}^n|z_{k+1}-z_k|^2.
\end{equation}
The constant is independent of $n$.
\end{lemma}
This is the discrete Hardy inequality in difference form; see
\citet[Chapter~IX, Theorem~326]{HLP52}.

\begin{lemma}[Expectations of random analytic functions]\label{lem:analyticexpectation}
Let $D\subset\C$ be a domain, and let $F(z)$, $z\in D$, be a jointly
measurable random function whose sample functions are almost surely
holomorphic on $D$. Suppose that, for every compact $K\subset D$,
\[
 \sup_{z\in K}\E|F(z)|<\infty.
\]
Then $f(z)=\E[F(z)]$ is holomorphic on $D$.
\end{lemma}
\begin{proof}
Choose a closed disk $\{|z-b|\le r\}\subset D$. For $|z-b|<r$,
Cauchy's integral formula and Fubini's theorem give
\[
 f(z)=\E\!\left[\frac1{2\pi\ii}\oint_{|w-b|=r}
                         \frac{F(w)}{w-z}\,dw\right]
     =\frac1{2\pi\ii}\oint_{|w-b|=r}
                         \frac{f(w)}{w-z}\,dw.
\]
The bound on $\E|F(w)|$ along the circle justifies Fubini's theorem. On every
smaller disk the denominator stays bounded away from zero, so the last
integral can be differentiated with respect to $z$. Thus $f$ is holomorphic
inside the disk. Such disks cover $D$.
\end{proof}

\begin{lemma}[Vitali's theorem]\label{lem:vitali}
Let $f_n$ be holomorphic on a connected open set $D\subset\C$ and locally
bounded there, uniformly in $n$. If $f_n$ converges on a set with an accumulation
point in $D$, then it converges locally uniformly on $D$ to a holomorphic function.
\end{lemma}
See \citet[Chapter~VII]{Conway78}.

\begin{lemma}[Two-constants theorem on a slit domain]\label{lem:twoconstants}
Let $D$ be a bounded smooth domain and let $J$ be a nondegenerate closed line
segment in $D$ such that $D\setminus J$ is connected. Let
$\omega(z)$ be the harmonic measure of $J$ in $D\setminus J$. If $f$ is
holomorphic on a neighborhood of $\overline D$, and
$|f|\le b_0$ on $J$, $|f|\le b_1$ on $\partial D$, where $b_0,b_1>0$, then
\begin{equation}\label{eq:twoconstants}
 |f(z)|\le b_0^{\omega(z)}b_1^{1-\omega(z)},\qquad z\in D\setminus J.
\end{equation}
Moreover $\omega$ is positive in $D\setminus J$ and extends continuously as $1$
on $J$. With this extension, $\inf_K\omega>0$ for every compact $K\subset D$.
\end{lemma}
This is the two-constants theorem in its harmonic-measure form, applied to
$\log|f|$; see \citet[Sections~4.2--4.3]{Ransford95}. The smooth outer
boundary and the line-segment slit are regular for the Dirichlet problem,
so the stated boundary continuity and positivity of $\omega$ apply.

\section{Fourier moments and the random walk}\label{sec:walk}
For real $\beta,\lambda>0$ and $t\le0$, write
$q_k(\beta,\lambda,t)=\E e^{\ii k\alpha_{\lambda,\beta}(t)}$.
Time translation of the phase SDE \eqref{eq:phase} gives
\begin{equation}\label{eq:momentscaling}
 q_k(\beta,\lambda,t)=q_k(\beta,\lambda e^{\beta t/4}).
\end{equation}
It\^o's formula applied to the angular exponential gives
\begin{equation}\label{eq:timemoments}
 \partial_t q_n(\beta,\lambda,t)
 =\frac{n^2}{2}(q_{n+1}+q_{n-1}-2q_n)(\beta,\lambda,t)
 +\ii n\lambda\frac\beta4e^{\beta t/4}q_n(\beta,\lambda,t).
\end{equation}
The noise integrand is bounded, so its stochastic integral has mean zero
on each finite time interval. The right-hand side is continuous by
bounded convergence.
By \eqref{eq:momentscaling}, equation~\eqref{eq:timemoments} gives
differentiability in $\lambda>0$ and \eqref{eq:momentode}.
The same scaling identity and the entrance condition imply
$q_k(\beta,\lambda)\to1$ as $\lambda\downarrow0$, by bounded convergence.

\begin{proof}[Proof of Proposition~\ref{prop:walk}]
The embedded simple random walk hits zero after finitely many jumps
almost surely \citep[Chapter~4]{Durrett19}. Thus $S$ is nonexplosive 
and $A_\beta$ is almost surely finite.
For $T>0$, applying Lemma~\ref{lem:fk} on $[-T,0]$ and using
\eqref{eq:momentscaling} gives
\begin{equation}\label{eq:twotime}
 q_n(\beta,\lambda)=\E_n\left[
  \exp\left(\ii\lambda\frac\beta4
       \int_0^T e^{-\beta t/4}S_t\,dt\right)
  q_{S_T}(\beta,\lambda e^{-\beta T/4})\right].
\end{equation}
For $T\ge\tau_0$, the terminal factor is $q_0=1$ and the exponential
is $e^{\ii\lambda A_\beta}$. Since $\tau_0<\infty$ almost surely and
the integrand has modulus at most one, bounded convergence as
$T\to\infty$ gives \eqref{eq:area}.
\end{proof}

The same Feynman--Kac argument gives uniqueness for the moment system
among uniformly bounded solutions.

\begin{corollary}[Uniqueness of bounded solutions]\label{prop:uniqueness}
For every real $\beta>0$, the system \eqref{eq:momentode}, with $q_0\equiv1$,
has a unique solution $(q_k)_{k\in\mathbb Z}$ of continuously differentiable
functions on $(0,\infty)$ satisfying
\[
 \sup_{k\in\mathbb Z,\,\lambda>0}|q_k(\lambda)|<\infty.
\]
\end{corollary}
\begin{proof}
The phase Fourier coefficients provide such a solution. Let
$r_k=q_k-\widetilde q_k$ be the difference of two solutions, so that
$r_0=0$ and $|r_k(\lambda)|\le c$ for some finite $c$.
For fixed $\lambda,T>0$, the functions
$t\mapsto r_k(\lambda e^{\beta t/4})$ satisfy \eqref{eq:timemoments}
on $[-T,0]$ by \eqref{eq:momentode}. Lemma~\ref{lem:fk} therefore gives
\[
 r_k(\lambda)=\E_k\!\left[
 \exp\!\left(\ii\lambda\frac\beta4
       \int_0^T e^{-\beta t/4}S_t\,dt\right)
 r_{S_T}(\lambda e^{-\beta T/4})\right],\qquad k\ge1.
\]
Since $r_0=0$ and $\tau_0<\infty$ almost surely,
\[
 |r_k(\lambda)|\le c\,\Pp_k(\tau_0>T)\longrightarrow0
 \qquad\text{as }T\to\infty.
\]
Applying the same argument to $\overline{r_{-k}}$ proves uniqueness
for negative indices.
\end{proof}

\section{Real Fourier decay from the walk}\label{sec:real}
We use the representation \eqref{eq:area} to prove decay of the phase
moments on the real parameter intervals needed for analytic continuation.
\begin{proposition}[Real Fourier decay]\label{prop:realdecay}
Let $J,K\subset(0,\infty)$ be compact sets for $\lambda$ and $\beta$,
respectively. There are $a,b>0$, depending only on $J,K$, such that
\begin{equation}\label{eq:realuniform}
 \sup_{\lambda\in J,\,\beta\in K}|q_k(\beta,\lambda)|
 \le b e^{-a|k|^{1/3}},\qquad k\in\mathbb Z.
\end{equation}
\end{proposition}
Consider the walk $S$ from Proposition \ref{prop:walk}. At height $k$, the walk fluctuates by order $k\sqrt h$ over a short time $h$,
producing area fluctuations of order $kh^{3/2}$. Taking $h=k^{-2/3}$
makes these fluctuations of order one. This motivates the study of the process on this time scale.

\begin{lemma}[Short-time path limit]\label{p:clt}
Put $h=k^{-2/3}$. If $k/2\le j\le 2k$ and $j/k\to x$ as
$k\to\infty$, then under $\Pp_j$,
\[
 \big(k^{-2/3}(S_{hu}-j):0\le u\le1\big)
 \ \Longrightarrow\ (xB_u:0\le u\le1),
\]
where $B$ is standard Brownian motion and convergence is in
$D([0,1])$ with its Skorokhod topology.
\end{lemma}
\begin{proof}
Let $Y$ be a simple symmetric continuous-time walk started at zero,
with rate $1/2$ to each neighbor. Changing its holding times realizes
the walk started at $j$ as
\[
 S_t=j+Y_{\int_0^t S_s^2\,ds},
\]
see \citet[Chapter~6, Section~1]{EK86} for the general
random time-change construction.
Set
\[
 W_k(v)=k^{-2/3}Y_{k^{4/3}v},\qquad
 c_k(u)=\int_0^u\left(\frac{S_{hv}}k\right)^2\,dv.
\]
Donsker's theorem \citep[Theorem~8.1.4]{Durrett19}, together with the
law of large numbers for the Poisson jump count, gives $W_k\Rightarrow B$
on every bounded time interval.

Let $\zeta_k$ be the exit time of $S$ from $[k/4,3k]$. Before exit,
$S_s/k\le3$, so the rescaled clock runs at rate at most $9$ and
\[
 \sup_{u\le1}\frac{|S_{hu\wedge\zeta_k}-j|}{k}
 \le k^{-1/3}\sup_{v\le9}|W_k(v)|\xrightarrow{\Pp}0.
\]
Since exit requires displacement at least $k/4$,
$\Pp_j(\zeta_k\le h)\to0$. Thus $S_{hu}/k\to x$ uniformly in
probability on $[0,1]$, and consequently $c_k(u)\to x^2u$ uniformly
in probability. The clock limit is deterministic, so joint convergence
and continuity of composition at continuous paths give
\[
 k^{-2/3}(S_{h\cdot}-j)=W_k\circ c_k
 \ \Longrightarrow\ B_{x^2\cdot}\stackrel d=xB.
\qedhere \]
\end{proof}

\begin{lemma}[Area before hitting zero]\label{lem:tube}
Let $D_n\Rightarrow \sigma B$ in $D([0,1])$ for some $\sigma>0$.
Suppose that almost surely $D_n$ cannot change sign without visiting
zero. Let $\xi_n$ be its first zero in $[1/2,1]$, or $1$ if there
is none. For every nonempty open interval $I\subset\mathbb R$,
\[
 \liminf_{n\to\infty}\Pp\left(\xi_n<1,\quad
                 \int_0^{\xi_n}D_n(u)\,du\in I\right)>0 .
\]
\end{lemma}
\begin{proof}
Choose $a$ and $r>0$ with $[a-r,a+r]\subset I$. On $[0,1/2]$ let $d$
be the triangular path from zero to $4a$ at time $1/4$ and back to
zero; its integral is $a$. On $[1/2,1]$ let $d$ run linearly from
zero to $r$ at time $3/4$ and then to $-r$ at time $1$.

On the tube $\|D_n-d\|_\infty<r/2$, the values at $3/4$ and $1$
have opposite signs, so $\xi_n<1$. Since $|d|\le r$ on the second half,
\[
 \left|\int_0^{\xi_n}D_n(u)\,du-a\right|
 \le \|D_n-d\|_\infty+\int_{1/2}^1 |d(u)|\,du < r .
\]
The Brownian limit has positive probability of lying in this tube:
$d$ is polygonal and starts at zero, so this follows from
Cameron--Martin and the positive probability of a ball about zero.
The map $f\mapsto\|f-d\|_\infty$ is continuous at continuous paths
in the Skorokhod topology. The continuous mapping theorem and the
open-set bound thus give the same positive lower bound for the
$\liminf$ of the tube probabilities.
\end{proof}

If the walk $S$ starts near a large $k$, then over time $h=k^{-2/3}$ its
area fluctuations are of order one. The next lemma turns these
fluctuations into a contraction of the weighted transition operator
in the supremum norm, by coupling paths with the same endpoint but
different accumulated phases.

\begin{lemma}[Uniform short-time contraction]\label{lem:contraction}
Fix $0<b_0\le b_1<\infty$ and $\ell_0\ge0$. Put
$h=k^{-2/3}$, $I_k=[k/2,2k]\cap\mathbb Z$, and
$\sigma_k=\inf\{t\ge0:S_t\notin I_k\}$.
For a real function $b$ on $[0,h]$ satisfying
$b_0\le b(t)\le b_1$ and $|b(t)-b(r)|\le\ell_0|t-r|$, define
\[
 (T_k^bf)(j)=\E_j\left[e^{\ii\int_0^h b(t)S_t\,dt}
                 f(S_h)\mathbf1_{\{\sigma_k>h\}}\right],\qquad j\in I_k.
\]
There are $\delta>0$ and $k_0$, depending only on $b_0,b_1,\ell_0$,
such that $\|T_k^b\|_{\infty\to\infty}\le1-\delta$ for $k\ge k_0$.
\end{lemma}
\begin{proof}
Start two independent walks $S_{1,t},S_{2,t}$ at $j\in I_k$.
Form coupled walks $\widehat S_{1,t},\widehat S_{2,t}$ by letting
the second follow the first after their first meeting at or after $h/2$.
On the independent paths put
\[
 D_k(u)=b(hu)k^{-2/3}(S_{1,hu}-S_{2,hu}),\qquad 0\le u\le1,
\]
and let $\xi_k$ be its first zero in $[1/2,1]$, or $1$ if none exists.
The coupled phase difference is $\int_0^{\xi_k}D_k(u)\,du$.
Uniformly over $j$ and $b$, there is $p>0$ such that
\[
 \Pp\left(\xi_k<1,\quad
  2\pi/3<\int_0^{\xi_k}D_k(u)\,du<4\pi/3\right)\ge p
\]
for all sufficiently large $k$. Indeed, otherwise take a sequence with
$k\to\infty$ and success probabilities tending to zero. Pass to a
subsequence with $j/k\to x\in[1/2,2]$ and $b(0)\to a\in[b_0,b_1]$.
Since $\sup_u|b(hu)-b(0)|\le\ell_0h$, Lemma~\ref{p:clt} and independence give
$D_k\Rightarrow a\sqrt2\,xB$. The difference cannot change sign without
visiting zero, since the independent walks almost surely never jump
simultaneously and $b$ is positive. Lemma~\ref{lem:tube} contradicts
the choice of sequence.

Let $G$ denote the event in the preceding display, so that $\Pp(G)\ge p$,
and let $\widehat\sigma_\ell$ be the exit time of
$\widehat S_{\ell,\cdot}$ from $I_k$. For $\|f\|_\infty\le1$,
extending $f$ by zero outside $I_k$, set
\[
 V_\ell=\exp\!\left(\ii\int_0^h b(t)\widehat S_{\ell,t}\,dt\right)
 f(\widehat S_{\ell,h})\mathbf1_{\{\widehat\sigma_\ell>h\}},
 \qquad \ell=1,2.
\]
The marginal laws give $\E V_1=\E V_2=(T_k^bf)(j)$, and $|V_\ell|\le1$.
On $G$, the walks agree from time $h\xi_k$ onward. If both survive
until $h$, they have the same endpoint and phase difference
$\Theta=\int_0^{\xi_k}D_k(u)\,du\in(2\pi/3,4\pi/3)$. Hence
\[
 |V_1+V_2|=|f(\widehat S_{1,h})|\,|1+e^{\ii\Theta}|
 \le2|\cos(\Theta/2)|\le1.
\]
If either walk is killed, at most one contribution is nonzero, so again
$|V_1+V_2|\le1$. Thus this bound holds throughout $G$, while
$|V_1+V_2|\le2$ always. Consequently,
\[
 |(T_k^bf)(j)|\le\tfrac12\E|V_1+V_2|
 \le\tfrac12\Pp(G)+\Pp(G^c)\le1-\frac p2.
\]
Taking the supremum over $j\in I_k$ and $\|f\|_\infty\le1$ proves
the assertion with $\delta=p/2$.
\end{proof}

\begin{lemma}[Exponential exit bound]\label{lem:escape}
For the walk started at $k$, every $T>0$, and $0<\theta\le1$,
\begin{equation}\label{eq:escape}
 \Pp_k(\sigma_k\le T)
 \le2\exp\left(-\frac{\theta k}{2}+4k^2T\theta^2\right).
\end{equation}
\end{lemma}
\begin{proof}
Put $f(i)=\cosh(\theta(i-k))$. For $i\in I_k$, the generator satisfies
\[
 (\mathcal Gf)(i)=i^2(\cosh\theta-1)f(i)
 \le4k^2\theta^2 f(i),
\]
since $\cosh\theta-1\le\theta^2$ for $0\le\theta\le1$.
Dynkin's formula stopped at $\sigma_k$ and Gronwall's inequality give
\[
 \E_k f(S_{T\wedge\sigma_k})\le e^{4k^2\theta^2T}.
\]
On $\{\sigma_k\le T\}$, the displacement at exit is at least $k/2$,
so $f(S_{T\wedge\sigma_k})\ge\cosh(\theta k/2)$. Therefore
\[
 \Pp_k(\sigma_k\le T)
 \le\frac{e^{4k^2\theta^2T}}{\cosh(\theta k/2)}
 \le2\exp\left(-\frac{\theta k}{2}+4k^2T\theta^2\right).
\qedhere \]
\end{proof}

\begin{proof}[Proof of Proposition~\ref{prop:realdecay}]
Since $q_0=1$ and $q_{-k}=\overline{q_k}$ for real $\beta$, it suffices
to take $k\ge1$. Fix $\lambda\in J$ and $\beta\in K$, and put
\[
 h=k^{-2/3},\qquad m=\lfloor k^{1/3}\rfloor,\qquad T=mh\le k^{-1/3}.
\]
Set $b(t)=\lambda(\beta/4)e^{-\beta t/4}$ and
$f(j)=q_j(\beta,\lambda e^{-\beta T/4})$. On $[0,1]$, these $b$ have common
positive lower and upper bounds and a common Lipschitz constant,
uniformly over $\lambda\in J$ and $\beta\in K$. Also $\|f\|_\infty\le1$.
By \eqref{eq:twotime},
\[
 q_k(\beta,\lambda)=\E_k\left[
 e^{\ii\int_0^T b(t)S_tdt}f(S_T)\right].
\]
For $r=0,\ldots,m-1$, define $b^{[r]}(t)=b(rh+t)$ on $[0,h]$.
Each $b^{[r]}$ satisfies the same bounds in Lemma~\ref{lem:contraction}.
Restricting to paths staying in $I_k$ until $T$, the Markov property gives
\[
 \E_k\left[e^{\ii\int_0^T b(t)S_tdt}f(S_T)
                  \mathbf1_{\{\sigma_k>T\}}\right]
 =(T_k^{b^{[0]}}T_k^{b^{[1]}}\cdots T_k^{b^{[m-1]}}f)(k).
\]
For $k\ge k_0$, where $k_0$ depends only on $J,K$,
Lemma~\ref{lem:contraction} therefore yields
\[
 |q_k(\beta,\lambda)|\le(1-\delta)^m+\Pp_k(\sigma_k\le T).
\]
Apply Lemma~\ref{lem:escape} with $\theta=\varepsilon k^{-2/3}$,
where $0<\varepsilon<1/8$ is fixed. Since $T\le k^{-1/3}$,
\[
 \Pp_k(\sigma_k\le T)
 \le2\exp\left[-(\varepsilon/2-4\varepsilon^2)k^{1/3}\right].
\]
Also $(1-\delta)^m\le c_1e^{-c_2k^{1/3}}$. Enlarging the prefactor to
cover the finitely many remaining $k$ proves \eqref{eq:realuniform}.
\end{proof}

\begin{corollary}[Smoothness of the phase density]\label{cor:realdensity}
For $\beta,\lambda>0$, the law of
$\alpha_{\lambda,\beta}(0)\bmod2\pi$ has a smooth density
$p_{\beta,\lambda}(\theta)$.
Its angular derivatives of every order are jointly continuous in
$(\beta,\lambda,\theta)$ and locally uniformly bounded.
\end{corollary}
\begin{proof}
By Proposition~\ref{prop:realdecay}, the Fourier series
\[
 p_{\beta,\lambda}(\theta)=\frac1{2\pi}
 \sum_{k\in\mathbb Z}q_k(\beta,\lambda)e^{-\ii k\theta}
\]
and all its termwise angular derivatives converge absolutely and
locally uniformly. By Fourier inversion
\citep[Chapter~I, \S7.1, Corollary, and \S7.4]{Katznelson04},
this defines a smooth density of the phase law.

The walk in \eqref{eq:area} has a law independent of $\beta$ and $\lambda$.
Each path is absorbed after finitely many jumps, so $A_\beta$ is almost
surely continuous in $\beta>0$. Thus, if
$(\beta_n,\lambda_n)\to(\beta,\lambda)$ in $(0,\infty)^2$, then
\[
 e^{\ii\lambda_n A_{\beta_n}}\longrightarrow e^{\ii\lambda A_\beta}
 \quad\text{almost surely},\qquad
 |e^{\ii\lambda_n A_{\beta_n}}|=1.
\]
Bounded convergence in \eqref{eq:area} gives joint continuity of each
$q_k$ in $(\beta,\lambda)$. Locally uniform convergence of every
termwise angular derivative series gives the claimed joint continuity
and local boundedness.
\end{proof}

\section{From boundary angles to correlation measures}\label{sec:pair}

\begin{proof}[Proof of Proposition~\ref{prop:anglelimit}]
Write $\alpha_\lambda=\alpha_{\lambda,\beta}(0)$. For a boundary angle $u$, put
\[
 \mu_u=\sum_{\lambda:\,\alpha_\lambda=u\bmod 2\pi}
       \frac{\delta_\lambda}{\partial_\lambda\alpha_\lambda}.
\]
For any continuous compactly supported $g$ on $\R$ and bounded
measurable $h$ on $\mathbb T$, changing variables
$u=\alpha_\lambda\bmod2\pi$ on each inverse branch gives
\[
 \int_{\mathbb T}h(u)\int_{\R}g(\lambda)\,\mu_u(d\lambda)\,du
 =\int_{\R}g(\lambda)h(\alpha_\lambda)\,d\lambda,
\]
since $du=\partial_\lambda\alpha_\lambda\,d\lambda$ cancels the spectral weight.
Thus, as measures on $\R$,
\[
 \int_{\mathbb T}h(u)\,\mu_u(d\lambda)\,du
 =h(\alpha_\lambda)\,d\lambda.
\]
Let $U_0=U$ be the uniform
angle in \eqref{eq:pointprocess}, and put $U_j=U+\varepsilon V_j\bmod 2\pi$,
where the $V_j$ are independent uniform variables on $[-1,1]$, independent
of $U$ and the phases. For $0<\varepsilon<\pi$ the angles have joint
density
\[\frac{1}{2\pi(2\varepsilon)^{k-1}}\qquad \mbox{on} \qquad
\Big\{\max_{1\le j\le k-1}\operatorname{dist}(u_j-u_0,2\pi\mathbb Z)<\varepsilon\Big\}
\]
and zero elsewhere.
By the substitution above, $\E\prod_{j=0}^{k-1}\mu_{U_j}(I_j)$ is therefore
the integral on the right of \eqref{eq:anglelimit}, including its prefactor,
with $\alpha_{\lambda_j}-\alpha_{\lambda_0}$ in place of $\alpha_{\lambda_j-\lambda_0}$.

The two vectors have the same joint law: subtract the equation for $\lambda_0$
in \eqref{eq:phase} from the equations for $\lambda_j$. The differences solve
the phase SDE \eqref{eq:phase} at $\lambda_j-\lambda_0$, all driven by the same rotated
Brownian motion $d\widetilde Z_t=e^{-\ii\alpha_{\lambda_0,\beta}(t)}dZ_t$.

Almost surely no endpoint of the $I_j$ has angle $U$, since $U$ is
independent and uniform. Only finitely many inverse branches meet each
interval near $U$, and their locations and weights vary continuously.
Thus $\mu_{U_j}(I_j)\to\mu_U(I_j)$ almost surely. Each $U_j$ is itself
uniform and independent of the phases. Given $n$ points in $I_j$,
the $r$th moment of $\mu_U(I_j)$ is at most $n^r\E w^r$ for $r\ge1$,
by convexity, where $w$ is one normalized weight. Proposition~\ref{prop:shooting}
and Lemma~\ref{lem:phasebounds} thus give moments of all orders;
H\"older's inequality bounds the products uniformly in $L^2$.
We may therefore pass to expectations. At the common angle $U$ every
tuple in the product consists of distinct points because the intervals
are disjoint. Given the points, their weights are independent with mean
one, so the limiting expectation is $\E\prod_{j=0}^{k-1}\Xi_\beta(I_j)$.
\end{proof}

\begin{proof}[Proof of Theorem~\ref{prop:pairphase}]
Take $k=2$ in Proposition~\ref{prop:anglelimit}, with
$\sup I_0<\inf I_1$. Writing $\lambda=\lambda_1-\lambda_0$, the quotient
in \eqref{eq:anglelimit} is
\[
 \frac1{2\varepsilon}\int_{-\varepsilon}^{\varepsilon}
                  p_{\beta,\lambda}(\theta)\,d\theta
 \longrightarrow p_{\beta,\lambda}(0)
\]
locally uniformly in $\lambda$, by Corollary~\ref{cor:realdensity}.
Thus $\rho_\beta(0,\lambda)=p_{\beta,\lambda}(0)/(2\pi)$ gives a
continuous version for $\lambda>0$. Translation invariance and symmetry
give the density off the diagonal, and the Fourier series in
Corollary~\ref{cor:realdensity} gives \eqref{eq:pairfourier}.
\end{proof}

\begin{remark}
For $N_\lambda=\Xi_\beta((0,\lambda])$, we note in passing that $\rho_\beta(0,\lambda)=\tfrac12\partial_\lambda^2\operatorname{Var}(N_\lambda)+1/(4\pi^2)$, where
\[
 \operatorname{Var}(N_\lambda)=\tfrac16\,\E_\nu[1-\cos(\lambda A_\beta)]
 +\tfrac{2}{\beta\pi^2}\,\E_1[\operatorname{Cin}(\lambda A_\beta)].
\]
Here $\operatorname{Cin}(x)=\int_0^x(1-\cos t)\,dt/t$, and $\nu(k)=6/(\pi^2k^2)$, $k\ge1$, is the initial law of the walk $S$ in the first term. The formula follows by conditioning on the phase in \eqref{eq:rounding},
using Fourier expansion, and applying It\^o's formula to $\alpha_\lambda^2$.
\end{remark}

\section{Holomorphic Fourier coefficients}\label{sec:holomorphic}
We now extend the Fourier moments to complex $\beta$ with $\Ree\beta>0$
using the representations \eqref{eq:area} and \eqref{eq:cutoff}.
We first construct a
holomorphic function for each cutoff, and then remove the cutoff.

\subsection{Holomorphy before removing the cutoff}
\begin{lemma}[Holomorphy and ODE for the cutoff moments]\label{lem:cutoffholomorphic}
For $1\le k\le n$, $\lambda>0$, and $\Ree\beta>0$, define
$q_{k,n}(\beta,\lambda)$ by \eqref{eq:cutoff}, with $A_\beta$ given by
the integral in \eqref{eq:area}. This function is holomorphic in $\beta$
and continuously differentiable in $\lambda$. It satisfies
\begin{equation}\label{eq:cutoffode}
 \lambda\partial_\lambda q_{k,n}
 =\ii k\lambda q_{k,n}
  +\frac{2k^2}{\beta}(q_{k+1,n}+q_{k-1,n}-2q_{k,n}),
 \quad q_{0,n}=1,\quad q_{n+1,n}=0.
\end{equation}
As $\lambda\downarrow0$,
\begin{equation}\label{eq:cutoffinitial}
 q_{k,n}(\beta,\lambda)\longrightarrow1-\frac{k}{n+1}.
\end{equation}
\end{lemma}
\begin{proof}
Write $E_n=\{\tau_0<\tau_{n+1}\}$. The walk and this event are independent
of $\beta$. On a retained path, $\tau_0<\infty$ and $S_v\le n$ before
absorption. The integral
\[
 A_\beta=\frac\beta4\int_0^\infty e^{-\beta v/4}S_v\,dv
\]
is supported on the finite interval $[0,\tau_0]$ and is entire in $\beta$.
This follows from Lemma~\ref{lem:analyticexpectation}
applied to normalized Lebesgue measure on $[0,\tau_0]$.
On the right half-plane it satisfies the path-independent bound
\begin{equation}\label{eq:areabound}
 |A_\beta|\le\frac{n|\beta|}{4}
           \int_0^\infty e^{-(\Ree\beta)v/4}\,dv
       =\frac{n|\beta|}{\Ree\beta}.
\end{equation}
Thus
\[
 |e^{\ii\lambda A_\beta}\mathbf1_{E_n}|
       \le\exp\!\left(\frac{n\lambda|\beta|}{\Ree\beta}\right).
\]
Lemma~\ref{lem:analyticexpectation} gives holomorphy of $q_{k,n}$ in
$\beta$, and dominated convergence gives continuity in $\lambda$.
The same bound and gambler's ruin give \eqref{eq:cutoffinitial}.

For real $\beta>0$, let $\widetilde S$ be the chain on
$\{0,\ldots,n+1\}$ obtained by absorbing $S$ at both endpoints:
$\widetilde S_v=S_{v\wedge\tau_0\wedge\tau_{n+1}}$.
Its generator $\mathcal G_n$ has zero rows at
$0$ and $n+1$. For $T>0$, define the finite-start cutoffs
$q^{(n)}_{-T,k}(t)$ on $[-T,0]$ by the finite system

\[
 \partial_t q^{(n)}_{-T,k}(t)=(\mathcal G_n q^{(n)}_{-T}(t,\cdot))(k)
              +\ii k\lambda\frac\beta4e^{\beta t/4}q^{(n)}_{-T,k}(t),
 \qquad q^{(n)}_{-T,k}(-T)=\mathbf1_{\{k\le n\}}.
\]
This finite system has a unique solution, which is bounded on $[-T,0]$,
with $q^{(n)}_{-T,0}\equiv1$ and $q^{(n)}_{-T,n+1}\equiv0$; the equations
hold at the endpoints as well as in the interior. Lemma~\ref{lem:fk},
with the time origin shifted to $-T$, gives
\[
 q^{(n)}_{-T,k}(t)=\E_k\!\left[
 \exp\!\left(\ii\lambda\frac\beta4 e^{\beta t/4}
            \int_0^{t+T}e^{-\beta v/4}\widetilde S_v\,dv\right)
 \mathbf1_{\{\widetilde S_{t+T}\ne n+1\}}\right].
\]
As $T\to\infty$, this converges to
\[
 q_{k,n}(\beta,\lambda e^{\beta t/4}).
\]
Indeed, the integrands agree once an endpoint has been reached, so
the difference of their expectations is bounded in modulus by
$2\Pp_k(\tau_0\wedge\tau_{n+1}>t+T)$.
Convergence is therefore uniform on compact time intervals.  Passing to the integral form of the finite ODE therefore
shows that $t\mapsto q_{k,n}(\beta,\lambda e^{\beta t/4})$ satisfies the same time-dependent system.
The change of variables $s=\lambda e^{\beta t/4}$ gives
\eqref{eq:cutoffode} for real $\beta>0$.

To extend the equation, fix $0<a<b$ and write its integral form:
\begin{equation}\label{eq:cutoffintegral}
 \begin{aligned}
 q_{k,n}(\beta,b)-q_{k,n}(\beta,a)
 &=\int_a^b\biggl[\ii kq_{k,n}(\beta,s)\\
 &\quad+\frac{2k^2}{\beta s}\bigl(q_{k+1,n}(\beta,s)
       +q_{k-1,n}(\beta,s)-2q_{k,n}(\beta,s)\bigr)\biggr]ds.
 \end{aligned}
\end{equation}
Both sides are holomorphic in $\beta$; the displayed path bound is uniform
for $s\in[a,b]$ on compact subsets of the half-plane and justifies the
integral by Lemma~\ref{lem:analyticexpectation}.
The identity theorem extends \eqref{eq:cutoffintegral} to $\Ree\beta>0$.
Its integrand is continuous in $s$, so the fundamental theorem of calculus
now proves differentiability and \eqref{eq:cutoffode}.
\end{proof}

\subsection{A uniform energy estimate}
The pathwise bound used above to prove holomorphy grows exponentially with
$n$. To remove the cutoff, we need a different estimate that is uniform in
$n$; this is the purpose of the following energy argument.
A related weighted energy argument is used for real-parameter uniqueness
in \citet[Proposition~16]{QV25}; here the estimate must be uniform in the
cutoff for complex $\beta$.

For $z=(z_1,\ldots,z_n)$ put
\[
 \|z\|_n^2=\sum_{k=1}^n\frac{|z_k|^2}{k^2},\qquad
 d_n(z)=|z_1|^2+\sum_{k=1}^{n-1}|z_{k+1}-z_k|^2+|z_n|^2.
\]
The endpoint values in $d_n$ are zero; the boundary value $q_{0,n}=1$ in the
ODE will instead appear as a forcing term.

\begin{lemma}[Energy bound]\label{lem:energy}
For every $n\ge1$, $\lambda>0$, and $\Ree\beta>0$,
\begin{equation}\label{eq:energybound}
 \sum_{k=1}^n\frac{|q_{k,n}(\beta,\lambda)|^2}{k^2}
 \le4\frac{|\beta|^2}{(\Ree\beta)^2}.
\end{equation}
\end{lemma}
\begin{proof}
Write $q_k=q_{k,n}(\beta,\lambda)$, $q=(q_1,\ldots,q_n)$, and $h(\lambda)=\|q(\lambda)\|_n^2$.
Multiplying \eqref{eq:cutoffode} by $\bar q_k/k^2$, summing, and taking real
parts gives, before summation by parts,
\begin{align*}
 \frac\lambda2 h'(\lambda)
 &=\Ree\left[\ii\lambda\sum_{k=1}^n\frac{|q_k|^2}{k}
       +\frac2\beta\sum_{k=1}^n
             \bar q_k(q_{k+1}+q_{k-1}-2q_k)\right]\\
 &=2\Ree\left[\beta^{-1}\sum_{k=1}^n
             \bar q_k(q_{k+1}+q_{k-1}-2q_k)\right].
\end{align*}
The imaginary potential contributes zero: its displayed sum is purely
imaginary. Summation by parts, with $q_0=1$ and $q_{n+1}=0$, now yields
\begin{equation}\label{eq:energyidentity}
 \lambda h'=-4a\,d_n(q)+4\Ree(\beta^{-1}\bar q_1),
 \qquad a=\Ree(\beta^{-1})>0.
\end{equation}
Put $b=|\beta^{-1}|$. Since $|q_1|\le\sqrt{d_n(q)}$,
\[
 4\Ree(\beta^{-1}\bar q_1)
 \le4b\sqrt{d_n(q)}
 \le2a\,d_n(q)+\frac{2b^2}{a}.
\]
Hardy's inequality (Lemma~\ref{lem:hardy}, equation~\eqref{eq:hardy})
gives $d_n(q)\ge h/4$, and hence
\begin{equation}\label{eq:energyinequality}
 \lambda h'\le-\frac a2 h+\frac{2b^2}{a}.
\end{equation}

By \eqref{eq:cutoffinitial},
\[
 \lim_{v\downarrow0}h(v)
 =\sum_{k=1}^n\frac{(1-k/(n+1))^2}{k^2}
 \le\sum_{k=1}^\infty\frac1{k^2}=\frac{\pi^2}{6}.
\]
For $0<\lambda_0<\lambda$, Gronwall's inequality applied in the
variable $\log\lambda$ (equivalently, multiplication by the integrating
factor $\lambda^{a/2}$) gives
\[
 h(\lambda)\le\left(\frac{\lambda_0}{\lambda}\right)^{a/2}h(\lambda_0)
 +\frac{4b^2}{a^2}\left[1-
             \left(\frac{\lambda_0}{\lambda}\right)^{a/2}\right].
\]
Let $\lambda_0\downarrow0$. The first term vanishes because
$h(\lambda_0)$ is bounded. Since $b/a=|\beta|/\Ree\beta$, this proves
\eqref{eq:energybound}.
\end{proof}

\subsection{Removing the cutoff}
\begin{proof}[Proof of Proposition~\ref{prop:holomorphic}]
For each fixed $k\ge1$, Lemma~\ref{lem:energy} gives
\[
 |q_{k,n}(\beta,\lambda)|\le2k\frac{|\beta|}{\Ree\beta},\qquad n\ge k.
\]
Thus the holomorphic cutoff transforms are locally bounded uniformly in $n$.
For real $\beta>0$, subtracting \eqref{eq:cutoff} from \eqref{eq:area}
and applying the gambler's ruin formula to the embedded simple random
walk \citep[Chapter~4]{Durrett19} gives
\[
 \begin{aligned}
 |q_k(\beta,\lambda)-q_{k,n}(\beta,\lambda)|
 &=\left|\E_k\!\left[e^{\ii\lambda A_\beta}
                 \mathbf1_{\{\tau_{n+1}<\tau_0\}}\right]\right|
 \le\Pp_k(\tau_{n+1}<\tau_0)=\frac{k}{n+1}.
 \end{aligned}
\]
This proves \eqref{eq:ruin} and hence
$q_{k,n}(\beta,\lambda)\to q_k(\beta,\lambda)$.
Vitali's theorem, Lemma~\ref{lem:vitali}, supplies a locally uniform
holomorphic limit on the whole half-plane. It is the unique extension of the
real Fourier coefficient.

For a fixed integer $m$, pass to the limit $n\to\infty$ in the first $m$ terms
of \eqref{eq:energybound}, and then let $m\to\infty$. We obtain
\begin{equation}\label{eq:fullenergy}
 \sum_{k=1}^\infty\frac{|q_k(\beta,\lambda)|^2}{k^2}
 \le4\frac{|\beta|^2}{(\Ree\beta)^2}.
\end{equation}
For $k\ge1$, define
\begin{equation}\label{eq:negativeextension}
 q_{-k}(\beta,\lambda)=\overline{q_k(\bar\beta,\lambda)},\qquad q_0(\beta,\lambda)=1.
\end{equation}
These are holomorphic and agree with the negative real Fourier moments.
Equation~\eqref{eq:fullenergy} proves~\eqref{eq:complexbound}. In particular,
\[
 |q_k(\beta,\lambda)|\le2|k|\frac{|\beta|}{\Ree\beta},\qquad k\ne0,
\]
uniformly in $\lambda>0$.
\end{proof}

\section{Analytic transfer and completion of the proof}\label{sec:transfer}
The complex-parameter energy bound and the real stretched-exponential
estimate combine through the two-constants theorem as follows.

\begin{lemma}[Transfer of stretched-exponential decay]\label{lem:transfer}
Let $f_k$, $k\ge1$, be holomorphic on $\Ree\beta>0$. Suppose that for some
$a\ge0$ and every compact $K$ in this half-plane,
\[
 \sup_{\beta\in K}|f_k(\beta)|\le c_K k^a.
\]
Suppose also that there is a nondegenerate compact interval
$J\subset(0,\infty)$ and constants $b,c>0$ such that
\[
 \sup_{\beta\in J}|f_k(\beta)|\le c e^{-b k^{1/3}}.
\]
Then for every compact $K\subset\{\Ree\beta>0\}$ there are constants
$b_K,c_K'>0$ such that
\begin{equation}\label{eq:transferred}
 \sup_{\beta\in K}|f_k(\beta)|\le c_K'e^{-b_K k^{1/3}}.
\end{equation}
The conclusion is uniform over an additional index whenever the constants in
the two hypotheses are uniform over that index.
\end{lemma}
\begin{proof}
Fix $K$. Choose a bounded smooth convex domain $D$ with
$K\cup J\subset D$ and $\overline D\subset\{\Ree\beta>0\}$.
The slit domain $D\setminus J$ is connected. Let $\omega$ be the harmonic
measure from Lemma~\ref{lem:twoconstants}, extended as $1$ on $J$, and put
$\omega_* =\inf_K\omega>0$.
The bounds on $\partial D$ and $J$, followed by \eqref{eq:twoconstants}, give
\[
 |f_k(\beta)|
 \le(c e^{-b k^{1/3}})^{\omega(\beta)}
       (c_{\overline D}k^a)^{1-\omega(\beta)}
 \le c_1 k^a e^{-b\omega_* k^{1/3}},\qquad \beta\in K.
\]
The same bound holds on $K\cap J$ directly. Absorbing the polynomial
factor gives \eqref{eq:transferred} with $b_K=b\omega_*/2$.
The proof preserves uniformity in an additional index.
\end{proof}

\begin{proof}[Proof of Proposition~\ref{prop:rapid}]
Fix a compact interval $L\subset(0,\infty)$ and apply
Lemma~\ref{lem:transfer} with $f_k(\beta)=q_k(\beta,\lambda)$ and $a=1$.
Proposition~\ref{prop:holomorphic} supplies the complex bound uniformly
in $\lambda>0$, and Proposition~\ref{prop:realdecay} supplies the real
stretched-exponential bound uniformly for $\lambda\in L$.
This proves the conclusion for positive indices. For negative indices,
use \eqref{eq:negativeextension} and the positive-index estimate on
$\{\bar\beta:\beta\in K\}$. Enlarging the prefactor covers $k=0$.
\end{proof}

\begin{proof}[Proof of Theorem~\ref{thm:main}]
Fix $\lambda>0$.
By Proposition~\ref{prop:rapid},
\begin{equation}\label{eq:finalcorrelation}
 \frac1{4\pi^2}\sum_{k\in\mathbb Z}q_k(\beta,\lambda)
\end{equation}
converges absolutely and locally uniformly on $\Ree\beta>0$.
Every summand is holomorphic by Proposition~\ref{prop:holomorphic}, so the sum
is holomorphic. For real $\beta>0$, Theorem~\ref{prop:pairphase}
identifies \eqref{eq:finalcorrelation} with $\rho_\beta(0,\lambda)$.
Uniqueness of the holomorphic extension follows
from the identity theorem.

To prove joint continuity, let $\lambda_n\to\lambda>0$.
By Proposition~\ref{prop:rapid}, the holomorphic functions
$\beta\mapsto\rho_\beta(0,\lambda_n)$ are locally bounded uniformly in $n$.
For real $\beta>0$, they converge to $\rho_\beta(0,\lambda)$ by
Theorem~\ref{prop:pairphase}. Vitali's theorem, Lemma~\ref{lem:vitali},
therefore gives this convergence locally uniformly on $\Ree\beta>0$.
Together with continuity in $\beta$ at fixed $\lambda$, this proves
joint continuity.
\end{proof}

\begin{corollary}[Holomorphic extension of the phase density]\label{cor:complexdensity}
For every fixed $\lambda>0$, the Fourier series
\[
 p_{\beta,\lambda}(\theta)=\frac1{2\pi}
        \sum_{k\in\mathbb Z}q_k(\beta,\lambda)e^{-\ii k\theta}
\]
extends the real phase density to $\Ree\beta>0$. It is holomorphic in $\beta$
and smooth in $\theta$; the series and every angular derivative converge
uniformly for $\beta$ in compact subsets of the half-plane and
$\theta\in\mathbb T$.
\end{corollary}
\begin{proof}
For every $j\ge0$, the bound in Proposition~\ref{prop:rapid} remains
summable after multiplication by $|k|^j$. Thus the differentiated series
converges locally uniformly, which also proves holomorphy.
\end{proof}


\section{\texorpdfstring
  {Fourier moments and phase density at $\beta=2$}
  {Fourier moments and phase density at beta=2}}
\label{sec:beta2decay}

We provide an explicit solution to \eqref{eq:momentode} at $\beta=2$. The solution shows that the stretched exponential
decay exponent $1/3$ in Proposition~\ref{prop:rapid} cannot be increased beyond $1/2$.
It also yields an explicit phase density, with a simple probabilistic
representation at integer multiples of $2\pi$.

For $a>0$, $b\in\mathbb R$, and $\Ree z>0$, Tricomi's confluent
hypergeometric function is defined by
\begin{equation}\label{eq:tricomi-definition}
U(a,b,z)=\frac{1}{\Gamma(a)}
\int_0^\infty e^{-zt}t^{a-1}(1+t)^{b-a-1}\,dt.
\end{equation}
We use its analytic continuation to $\mathbb C\setminus(-\infty,0]$,
with complex powers on the principal branch, and set $U(0,b,z)=1$.
See \citet[\href{https://dlmf.nist.gov/13.4.E4}{(13.4.4)}]{DLMF}.

\begin{proposition}[Fourier moments at $\beta=2$]
For $k\ge1$ and $\lambda>0$,
\begin{equation}\label{eq:explicit}
q_k(2,\lambda)=k!\bigl[U(k,1,-\ii\lambda)
        -k e^{\ii\lambda}U(k+1,1,-\ii\lambda)\bigr].
\end{equation}
For fixed
$\lambda>0$, as $k\to\infty$,
\begin{equation}\label{eq:decay}
|q_k(2,\lambda)|\sim
\begin{cases}
2\sqrt\pi\,\lambda^{-1/4}|\sin(\lambda/2)|\,k^{3/4}
 e^{-\sqrt{2\lambda k}}, & \lambda\notin2\pi\mathbb Z,\\[2pt]
\sqrt\pi\,\lambda^{1/4}k^{1/4}e^{-\sqrt{2\lambda k}},
 & \lambda\in2\pi\mathbb Z.
\end{cases}
\end{equation}
\end{proposition}

\begin{proof}
Let $r_k(\lambda)$ be the right-hand side of \eqref{eq:explicit}, with
$r_0=1$ and $r_{-k}=\overline{r_k}$. The recurrence and derivative
identities for $U$ in DLMF
\href{https://dlmf.nist.gov/13.3.E7}{(13.3.7)},
\href{https://dlmf.nist.gov/13.3.E10}{(13.3.10)}, and
\href{https://dlmf.nist.gov/13.3.E22}{(13.3.22)} give
\[
\lambda r_k'=\ii k\lambda r_k+
 k^2(r_{k+1}+r_{k-1}-2r_k),\qquad k\ge1.
\]
The negative-index equations follow by conjugation.

For $\Ree z>0$, \eqref{eq:tricomi-definition} and a shift in the
second integral give
\begin{align}
&k!\bigl[U(k,1,z)-k e^{-z}U(k+1,1,z)\bigr]
 =\int_0^\infty e^{-zt}h_k(t)\,dt,\label{eq:integral}\\
&h_k(t)=k\left[
 \frac{t^{k-1}}{(1+t)^k}
 -\mathbf1_{\{t>1\}}\frac{(t-1)^k}{t^{k+1}}
\right].\notag
\end{align}
For $t>1$ the ratio of the second term in brackets to the first is
$(1-t^{-2})^k$, so $h_k\ge0$.
Integrating up to $R>1$ and substituting $s=t-1$ in the second
term leaves a boundary term bounded by $k/R$; letting $R\to\infty$ gives
\[
\int_0^\infty h_k(t)\,dt
=k\int_0^\infty\frac{t^{k-1}}{(1+t)^{k+1}}\,dt=1.
\]
Letting $z=\varepsilon-\ii\lambda$ with $\varepsilon\downarrow0$
in \eqref{eq:integral}, dominated convergence gives
$r_k(\lambda)=\int_0^\infty e^{\ii\lambda t}h_k(t)\,dt$.
Hence $|r_k(\lambda)|\le1$, and Corollary~\ref{prop:uniqueness}
identifies $r_k$ with $q_k(2,\lambda)$.

The recurrence identity (DLMF \href{https://dlmf.nist.gov/13.3.E9}{(13.3.9)})
\[
U(k,1,z)-kU(k+1,1,z)=U(k,0,z)
\]
rewrites the left-hand side of \eqref{eq:integral} as
\[
k!\bigl[(1-e^{-z})U(k,1,z)+e^{-z}U(k,0,z)\bigr].
\]
On the principal branch, the large-parameter expansion of $U$ and
the large-argument expansion of the modified Bessel function
$K_\nu$ of the second kind in DLMF
\href{https://dlmf.nist.gov/13.8.E11}{(13.8.11)} and
\href{https://dlmf.nist.gov/10.40.E2}{(10.40.2)} give
\[
k!U(k,b,z)\sim
\sqrt\pi\,e^{z/2}z^{1/4-b/2}k^{1/4+b/2}e^{-2\sqrt{kz}},
\qquad b\in\{0,1\},
\]
locally uniformly for $z\ne0$ with $\Ree z\ge0$.
At $z=-\ii\lambda$ the $b=1$ term dominates unless $e^{\ii\lambda}=1$,
when it vanishes. This gives \eqref{eq:decay}.
\end{proof}

\begin{samepage}
The same Fourier coefficients give an explicit phase density. We use the
exponential integral
\[
 E_1(w)=\int_1^\infty \frac{e^{-wt}}{t}\,dt,
 \qquad \Ree w>0,
\]
continued to $\C\setminus(-\infty,0]$ on its principal branch; see
\citet[\href{https://dlmf.nist.gov/6.2.i}{\S6.2(i)}]{DLMF}.
\end{samepage}

\begin{corollary}[Phase density and pair correlation at $\beta=2$]\label{cor:beta2density}
For $\lambda>0$ and $0<\theta<2\pi$, put
\begin{equation}\label{eq:beta2w}
 w=-\frac{\ii\lambda}{1-e^{-\ii\theta}}
   =-\frac{\lambda}{2}\cot\frac\theta2-\frac{\ii\lambda}{2}.
\end{equation}
The density of $\alpha_{\lambda,2}(0)\bmod2\pi$ with respect to $d\theta$ is
\begin{equation}\label{eq:beta2density}
 p_{2,\lambda}(\theta)
 =\frac{1-\cos\lambda+
 \Ree\!\left(\bigl[(1+w)(e^{\ii\lambda}-1)-\ii\lambda\bigr]
                  e^wE_1(w)\right)}
 {2\pi(1-\cos\theta)}.
\end{equation}
Moreover, we recover the classical sine-kernel formula \citep{Dyson1962}
\[
\rho_2(0,\lambda)=\frac1{4\pi^2}
\left[1-\left(\frac{\sin(\lambda/2)}{\lambda/2}\right)^2\right].
\]
\end{corollary}

\begin{proof}
Put $\eta=e^{-\ii\theta}$. For $\Ree z>0$, summing the geometric series
in \eqref{eq:integral} and shifting the second integral gives
\[
 \sum_{k\ge1}\eta^k k!\bigl[U(k,1,z)-ke^{-z}U(k+1,1,z)\bigr]
 =\eta\int_0^\infty
 \frac{e^{-zt}\bigl[1+(1-e^{-z})t\bigr]}
      {\bigl[1+(1-\eta)t\bigr]^2}\,dt.
\]
The interchange is justified by $h_k\ge0$ and
$\sum_{k\ge1}h_k(t)=1+\min(t,1)$. Integrating by parts and using
\[
 \int_0^\infty\frac{e^{-zt}}{1+(1-\eta)t}\,dt
 =\frac{e^{z/(1-\eta)}E_1(z/(1-\eta))}{1-\eta},
\]
we may let $z\to-\ii\lambda$ by the locally uniform decay above to obtain
\[
 \sum_{k\ge1}q_k(2,\lambda)e^{-\ii k\theta}
 =\frac{e^{-\ii\theta}-e^{\ii\lambda}
 +\bigl[(1+w)(e^{\ii\lambda}-1)-\ii\lambda\bigr]e^wE_1(w)}
 {2(1-\cos\theta)}.
\]
Fourier inversion gives \eqref{eq:beta2density}, and the decay of the
Fourier coefficients gives continuity at $\theta=0$.

At angle zero, the identity $\sum_{k\ge1}h_k(t)=1+\min(t,1)$ gives,
for $\Ree z>0$,
\[
\int_0^\infty e^{-zt}\bigl(1+\min(t,1)\bigr)\,dt
=\frac1z+\frac{1-e^{-z}}{z^2}.
\]
The locally uniform decay in $k$ allows passage to $z=-\ii\lambda$,
giving 
\[
 \sum_{k\ge1}q_k(2,\lambda)
 =\frac{\ii}{\lambda}
   +\frac{e^{\ii\lambda}-1}{\lambda^2}.
\]
Substituting into the Fourier formula \eqref{eq:finalcorrelation}
\[
\rho_2(0,\lambda)=\frac1{4\pi^2}
  \left(1+2\Ree\sum_{k\ge1} q_k(2,\lambda)\right)
\]
gives the stated sine-kernel formula.
\end{proof}

\begin{remark}[The phase at integer multiples of $2\pi$]
At $\lambda=2\pi m$, $m\ge1$ an integer, the density
\eqref{eq:beta2density} simplifies to the positive real integral
\begin{equation}\label{eq:beta2integerdensity}
 p_{2,2\pi m}(\theta)=\frac m2\int_0^\infty
 \frac{e^{-\pi m t}}
      {1-t\sin\theta+t^2\sin^2(\theta/2)}\,dt,
 \qquad 0<\theta<2\pi.
\end{equation}
Equivalently,
\begin{equation}\label{eq:beta2cauchyexp}
 \alpha_{2\pi m,2}(0)\bmod2\pi
 \ \stackrel d=\ 2\operatorname{arccot}(X+Y) \bmod2\pi,
\end{equation}
where $X$ is standard Cauchy and $Y$ is independent exponential with rate
$\pi m$.

To see this, substitute $e^{\ii\lambda}=1$ in
\eqref{eq:beta2density} and use
\[
 e^wE_1(w)=\int_0^\infty\frac{e^{-s}}{s+w}\,ds,
 \qquad w\notin(-\infty,0],
\]
the analytic continuation of DLMF
\href{https://dlmf.nist.gov/6.7.E1}{(6.7.1)}.
Evaluating the real part and setting $s=\pi m t$ gives
\eqref{eq:beta2integerdensity}. Under $x=\cot(\theta/2)$, the density becomes the convolution of
the exponential density with rate $\pi m$ and the standard Cauchy
density, proving \eqref{eq:beta2cauchyexp}.

\end{remark}

\medskip
\noindent\textbf{AI use.} ChatGPT was used extensively to develop the random walk coupling,
the pair correlation identity, and the analytic continuation argument, and to
check calculations, locate references, and revise the manuscript. The authors
take full responsibility for the mathematics and references.

\medskip
\noindent\textbf{Acknowledgments.}
Benedek Valk\'o was partially supported by the National Science
Foundation award DMS-2246435. B\'alint Vir\'ag was partially supported by an NSERC
Discovery Grant.

\setlength{\bibsep}{3pt}

\par\addvspace{2\baselineskip}
\begingroup
\normalsize\raggedright
\noindent Yun Li. School of Mathematical Sciences, Nankai University, Tianjin, China.\\*
\href{mailto:liyun1723@hotmail.com}{\texttt{liyun1723@hotmail.com}}\par
\medskip
\noindent Yahui Qu. Department of Mathematics, University of Wisconsin--Madison,
\mbox{Madison, USA.}\\*
\href{mailto:yqu45@wisc.edu}{\texttt{yqu45@wisc.edu}}\par
\medskip
\noindent Benedek Valk\'o. Department of Mathematics, University of Wisconsin--Madison,
\mbox{Madison, USA.}\\*
\href{mailto:valko@math.wisc.edu}{\texttt{valko@math.wisc.edu}}\par
\medskip
\noindent B\'alint Vir\'ag. Department of Mathematics, University of Toronto, Toronto, Canada.\\*
\href{mailto:balint@math.toronto.edu}{\texttt{balint@math.toronto.edu}}\par
\endgroup

\begin{thebibliography}{99}
\bibitem[Assiotis and Najnudel(2026)]{AN26}
T.~Assiotis and J.~Najnudel (2026).
\newblock Moments of C$\beta$E field partition function,
$\operatorname{Sine}_\beta$ correlations and stochastic zeta.
\newblock \href{https://arxiv.org/abs/2602.08739}{arXiv:2602.08739}.

\bibitem[Conway(1978)]{Conway78}
J.~B.~Conway (1978).
\newblock \emph{Functions of One Complex Variable I}, second edition.
\newblock Graduate Texts in Mathematics \textbf{11}, Springer.

\bibitem[DLMF(2026)]{DLMF}
\emph{NIST Digital Library of Mathematical Functions} (2026).
\newblock Release 1.2.8 of 2026-09-15.
\newblock F.~W.~J.~Olver, A.~B.~Olde Daalhuis, D.~W.~Lozier,
B.~I.~Schneider, R.~F.~Boisvert, C.~W.~Clark, B.~R.~Miller,
B.~V.~Saunders, H.~S.~Cohl, and M.~A.~McClain, eds.
\newblock \url{https://dlmf.nist.gov/}.

\bibitem[Durrett(2019)]{Durrett19}
R.~Durrett (2019).
\newblock \emph{Probability: Theory and Examples}, fifth edition.
\newblock Cambridge University Press.

\bibitem[Dyson(1962)]{Dyson1962}
F.~J. Dyson  (1962).
\newblock Statistical theory of the energy levels of complex systems.~III.
\newblock \emph{J. Math. Phys.} \textbf{3}, 166--175.

\bibitem[Ethier and Kurtz(1986)]{EK86}
S.~N.~Ethier and T.~G.~Kurtz (1986).
\newblock \emph{Markov Processes: Characterization and Convergence}.
\newblock John Wiley \& Sons.

\bibitem[Fitzsimmons and Pitman(1999)]{FP99}
P.~J.~Fitzsimmons and J.~Pitman (1999).
\newblock Kac's moment formula and the Feynman--Kac formula for additive functionals of a Markov process.
\newblock \emph{Stochastic Processes and their Applications} \textbf{79}, 117--134.
\newblock \url{https://www.stat.berkeley.edu/~pitman/kac.pdf}.

\bibitem[Forrester(1992)]{Forrester92}
P.~J.~Forrester (1992).
\newblock Selberg correlation integrals and the $1/r^2$ quantum many-body system.
\newblock \emph{Nuclear Physics B} \textbf{388}(3), 671--699.

\bibitem[Forrester(1994)]{Forrester94}
P.~J.~Forrester (1994).
\newblock Addendum to ``Selberg correlation integrals and the $1/r^2$ quantum many-body system''.
\newblock \emph{Nuclear Physics B} \textbf{416}(1), 377--385.

\bibitem[Forrester(2010)]{Forrester10}
P.~J.~Forrester (2010).
\newblock \emph{Log-Gases and Random Matrices}.
\newblock London Mathematical Society Monographs \textbf{34},
Princeton University Press.

\bibitem[Forrester, Jancovici and McAnally(2001)]{FJM01}
P.~J.~Forrester, B.~Jancovici and D.~S.~McAnally (2001).
\newblock Analytic properties of the structure function for the one-dimensional one-component log-gas.
\newblock \emph{Journal of Statistical Physics} \textbf{102}(3--4), 737--780.
\newblock \href{https://doi.org/10.1023/A:1004846818738}{doi:10.1023/A:1004846818738}.

\bibitem[Ha(1995)]{Ha95}
Z.~N.~C.~Ha (1995).
\newblock Fractional statistics in one dimension: view from an exactly solvable model.
\newblock \emph{Nuclear Physics B} \textbf{435}(3), 604--636.
\newblock \href{https://doi.org/10.1016/0550-3213(94)00537-O}{doi:10.1016/0550-3213(94)00537-O}.

\bibitem[Hardy, Littlewood and P\'olya(1952)]{HLP52}
G.~H.~Hardy, J.~E.~Littlewood and G.~P\'olya (1952).
\newblock \emph{Inequalities}, second edition.
\newblock Cambridge University Press.

\bibitem[Jancovici(1981)]{Jancovici81}
B.~Jancovici (1981).
\newblock Exact results for the two-dimensional one-component plasma.
\newblock \emph{Physical Review Letters} \textbf{46}(6), 386--388.
\newblock \href{https://doi.org/10.1103/PhysRevLett.46.386}{doi:10.1103/PhysRevLett.46.386}.

\bibitem[Katznelson(2004)]{Katznelson04}
Y.~Katznelson (2004).
\newblock \emph{An Introduction to Harmonic Analysis}, third edition.
\newblock Cambridge University Press.

\bibitem[Killip and Stoiciu(2009)]{KS09}
R.~Killip and M.~Stoiciu (2009).
\newblock Eigenvalue statistics for CMV matrices: From Poisson to clock via random matrix ensembles.
\newblock \emph{Duke Mathematical Journal} \textbf{146}, 361--399.

\bibitem[Li and Valk\'o(2022)]{LV_HP}
Y.~Li and B.~Valk\'o (2022).
\newblock Operator level limit of the circular Jacobi $\beta$-ensemble.
\newblock \emph{Random Matrices: Theory and Applications} \textbf{11}(4), 2250043.
\newblock \href{https://doi.org/10.1142/S2010326322500435}{doi:10.1142/S2010326322500435}.

\bibitem[Nakano(2014)]{Nakano14}
F.~Nakano (2014).
\newblock Level statistics for one-dimensional Schr\"odinger operators and Gaussian beta ensemble.
\newblock \emph{Journal of Statistical Physics} \textbf{156}, 66--93.
\newblock \href{https://doi.org/10.1007/s10955-014-0987-x}{doi:10.1007/s10955-014-0987-x}.

\bibitem[Okounkov(1997)]{Okounkov97}
A.~Okounkov (1997).
\newblock On $n$-point correlations in the log-gas at rational temperature.
\newblock \href{https://arxiv.org/abs/hep-th/9702001}{arXiv:hep-th/9702001}.


\bibitem[Qu and Valk\'o(2025)]{QV25}
Y.~Qu and B.~Valk\'o (2025).
\newblock On the pair correlation function of the $\operatorname{Sine}_\beta$ process.
\newblock \href{https://arxiv.org/abs/2509.15446v2}{arXiv:2509.15446v2}, revised October 2026.

\bibitem[Ransford(1995)]{Ransford95}
T.~Ransford (1995).
\newblock \emph{Potential Theory in the Complex Plane}.
\newblock London Mathematical Society Student Texts \textbf{28}, Cambridge University Press.

\bibitem[Valk\'o and Vir\'ag(2009)]{VV09}
B.~Valk\'o and B.~Vir\'ag (2009).
\newblock Continuum limits of random matrices and the Brownian carousel.
\newblock \emph{Inventiones Mathematicae} \textbf{177}, 463--508.

\bibitem[Valk\'o and Vir\'ag(2017)]{VV17}
B.~Valk\'o and B.~Vir\'ag (2017).
\newblock The $\operatorname{Sine}_\beta$ operator.
\newblock \emph{Inventiones Mathematicae} \textbf{209}, 275--327.

\bibitem[Valk\'o and Vir\'ag(2022)]{VV22}
B.~Valk\'o and B.~Vir\'ag (2022).
\newblock The many faces of the stochastic zeta function.
\newblock \emph{Geometric and Functional Analysis} \textbf{32}, 1160--1231.
\newblock \href{https://doi.org/10.1007/s00039-022-00613-8}{doi:10.1007/s00039-022-00613-8}.

\bibitem[Valk\'o and Vir\'ag(2023)]{VV23}
B.~Valk\'o and B.~Vir\'ag (2023).
\newblock Palm measures for Dirac operators and the $\operatorname{Sine}_\beta$ process.
\newblock \emph{Stochastic Processes and their Applications} \textbf{163}, 106--135.
\newblock \href{https://doi.org/10.1016/j.spa.2023.05.011}{doi:10.1016/j.spa.2023.05.011}.
\end{thebibliography}
\end{document}